\documentclass{amsart}

\usepackage{xy,amssymb,diagrams}
\usepackage[hyperindex,pdfmark]{hyperref}

\xyoption{poly}
\xyoption{arc}
\xyoption{knot}
\xyoption{curve}
\xyoption{arrow}
\xyoption{all}

\makeindex

\newcommand{\een}{{\xy 0;/r.08pc/:
(3,-3); (5,-3) **@{-}; (5,5) **@{-}; (0,2) **@{-};
(3,-3); (3,4) **@{-}
\endxy}}

\newcommand{\dummy[1]}{{#1}}
\newcommand{\C}{\mathbb{C}}

\newcommand{\Q}{\mathbb{Q}}

\newcommand{\N}{\mathbb{N}}
\newcommand{\Z}{\mathbb{Z}}

\newarrow{Implies} ===={=>}

\newcommand{\wis}[1]{{\text{\em \usefont{OT1}{cmtt}{m}{n} #1}}}

\newtheorem{theorem}{Theorem}[section]

\newtheorem{proposition}[theorem]{Proposition}

\theoremstyle{definition}

\theoremstyle{remark}

\begin{document}

\sloppy

\title{Noncommutative smoothness and coadjoint orbits}

\author{Lieven Le Bruyn}
\address{Universiteit Antwerpen (UIA) \\ B-2610 Antwerp (Belgium)}
\email{lebruyn@uia.ua.ac.be}
\urladdr{http://win-www.uia.ac.be/u/lebruyn/}

\maketitle

\begin{center} {\it For Claudio Procesi on his 60th birthday.}
\end{center}

\begin{abstract}
In \cite{BocklandtLB:2000} R. Bocklandt and the author proved that certain quotient varieties
of representations of deformed preprojective algebras are coadjoint orbits for the necklace
Lie algebra $\mathfrak{N}_Q$ of the corresponding quiver $\vec{Q}$. A conjectural ringtheoretical explanation of these results was given in terms of
noncommutative smoothness in the sense of C. Procesi \cite{Procesi:1987}. In this paper we 
prove these conjectures. The main tool in the proof is the \'etale local description due to
W. Crawley-Boevey \cite{Crawley:2001}. Along the way we determine the smooth locus of the
Marsden-Weinstein reductions for quiver representations.
\end{abstract}

\section{Introduction.}

In \cite{BerestWilson:2000} Yu. Berest and G. Wilson asked whether the Calogero-Moser phase
space is a coadjoint orbit for a central extension of the automorphism group of the Weyl algebra.
This is indeed the case as was first proved by V. Ginzburg \cite{Ginzburg:1999} and subsequently generalized independently
by V. Ginzburg \cite{Ginzburg:2000} and R. Bocklandt and the author \cite{BocklandtLB:2000} to
certain quiver-varieties. Both proofs use noncommutative symplectic geometry as outlined by
M. Kontsevich \cite{Kontsevich:1993} in an essential way.

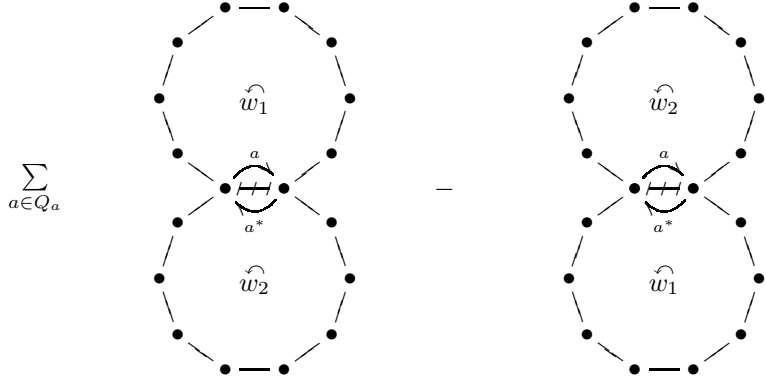
\begin{figure}
\[
\begin{xy}/r3pc/:
{\xypolygon10{~*{\bullet}~>>{}}},
"1" *+{\txt{\tiny $u$}},"8"="a","10"="a1","9"="a2","0"="c1",
"c1"+(0,-1.9),
{\xypolygon10{~*{\bullet}~>>{}}},
"5" *+{\txt{\tiny $v$}},"4"="b","6"="b1","5"="b2","0"="c2",
"a";"a2" **@{/};
"a" +(-2,0) *+{\txt{$\underset{a \in Q_a}{\sum}$}};
"c1" *+{\txt{$\overset{\curvearrowleft}{w_1}$}};
"c2" *+{\txt{$\overset{\curvearrowleft}{w_2}$}};
\POS"a" \ar@/^2ex/^{\txt{\tiny{$a$}}} "a2"
\POS"a2" \ar@/^2ex/^{\txt{\tiny{$a^*$}}} "a"
\end{xy}~\qquad
  \begin{xy}/r3pc/:
{\xypolygon10{~*{\bullet}~>>{}}},
"1" *+{\txt{\tiny $u$}},"8"="a","10"="a1","9"="a2","0"="c1",
"c1"+(0,-1.9),
{\xypolygon10{~*{\bullet}~>>{}}},
"5" *+{\txt{\tiny $v$}},"4"="b","6"="b1","5"="b2","0"="c2",
"a";"a2" **@{/};
"a"+(-2,0) *+{-};
"c1" *+{\txt{$\overset{\curvearrowleft}{w_2}$}};
"c2" *+{\txt{$\overset{\curvearrowleft}{w_1}$}};
\POS"a" \ar@/^2ex/^{\txt{\tiny{$a$}}} "a2"
\POS"a2" \ar@/^2ex/^{\txt{\tiny{$a^*$}}} "a"
\end{xy}
\]
\caption{Lie bracket $[ w_1,w_2 ]$ in $\mathfrak{N}_Q$.}
\label{bracket}
\end{figure}

Recall that a quiver $\vec{Q}$ is a finite directed graph on a set of vertices $Q_v = \{ v_1,\hdots,
v_k \}$ and having a finite set of arrows $Q_a = \{ a_1,\hdots,a_l \}$ where we allow both multiple
arrows between vertices and loops in vertices. The double quiver $\bar{Q}$ of the
quiver $\vec{Q}$ is the quiver obtained by adjoining to every arrow $a \in Q_a$ an arrow $a^*$ in the
opposite direction. Two oriented cycles in $\bar{Q}$ are equivalent if they are equal up to a cyclic
permutation of the arrow components. A {\em necklace word} $w$ for $\bar{Q}$ is an equivalence class of
oriented cycles in $\bar{Q}$. The {\em necklace Lie algebra} $\mathfrak{N}_Q$ of the quiver
$\vec{Q}$ has as basis the set of all necklace words $w$ for $\bar{Q}$ and with Lie bracket
$[w_1,w_2]$ determined by figure~\ref{bracket}.
That is, for every arrow $a \in Q_a$ we look for an occurrence of $a$ in
$w_1$ and of $a^*$ in
$w_2$. We then open up the necklaces by removing these factors and regluing
the open ends together
to form a new necklace word. We repeat this operation for {\em all}
occurrences of $a$ (in $w_1$)
and $a^*$ (in $w_2$). We then replace the roles of $a^*$ and $a$ and redo
this operation with a
minus sign. Finally, we add up all these obtained necklace words for all
arrows $a \in Q_a$.

The path algebra $\C \bar{Q}$ has as $\C$-basis the set of all oriented paths
$p = a_{i_u} \hdots a_{i_1}$ of length $u \geq 1$ together with the vertex-idempotents
$e_i$ considered as paths of length zero. Multiplication in $\C \bar{Q}$ is induced
by concatenation (on the left) of paths. Let $V = \C \times \hdots \times \C$ be the
$k$-dimensional semisimple subalgebra generated by the vertex-idempotents. In
\cite{BocklandtLB:2000} the noncommutative {\em relative} differential forms $\Omega^i_V \C \bar{Q}$
(introduced and studied by J. Cuntz and D. Quillen in \cite{CuntzQuillen:1995}) were used to
describe the noncommutative relative deRham (or Karoubi) complex
\[
\wis{dR}^0_V~\C Q \rTo^d \wis{dR}^1_V~\C Q \rTo^d \wis{dR}^2_V~\C Q \rTo^d
\hdots \]
where we define the vectorspace quotients dividing out the super-commutators
\[
\wis{dR}^n_V~\C Q = \dfrac{\Omega^n_V~\C Q}{\sum_{i=0}^n [~\Omega^i_V~\C Q,
\Omega^{n-i}_V~\C Q~]} \]
In particular, the noncommutative functions $\wis{dR}^0_V~\C \bar{Q}$ coincide with $\mathfrak{N}_Q$.
A noncommutative symplectic structure is defined on $\C \bar{Q}$ by the element
$\omega = \sum_{a \in Q_a} da^* da \in \wis{dR}^2_V~\C \bar{Q}$ and we have a noncommutative version
of the classical result in symplectic geometry relating the Lie algebra of functions to Hamiltonian
vectorfields : there is a central extension of Lie algebras
\[
0 \rTo V \rTo \mathfrak{N}_Q \rTo Der_{\omega}~\C \bar{Q} \rTo 0 \]
where $Der_{\omega}~\C \bar{Q}$ is the Lie algebra of {\it symplectic} derivations, that is,
 $\theta \in Der_V~\C \bar{Q}$ such that $L_{\theta} \omega = 0$ where $L_{\theta}$ is the degree
 preserving derivation on the relative differential forms determined by $L_{\theta}(a) = \theta(a)$ and
 $L_{\theta}(da) = d \theta(a)$, see \cite[Thm. 4.2]{BocklandtLB:2000}. As $Der_{\omega}~\C \bar{Q}$
 corresponds to the group of $V$-automorphisms of $\C \bar{Q}$ preserving the element
 $m = \sum_{a \in Q_a} [a,a^*]$ it is natural to consider for $\lambda = \sum_{i=1}^k \lambda_i e_i$
 with $\lambda_i \in \Q$ the {\em deformed preprojective algebra}
 \[
 \Pi_{\lambda}(\bar{Q}) = \dfrac{\C \bar{Q}}{(m-\lambda)} \]
 For a given dimension vector $\alpha = (a_1,\hdots,a_k) \in \N^k$ one defines the affine scheme
 $\wis{rep}_{\alpha}~\Pi_{\lambda}$ of $\alpha$-dimensional representations of $\Pi_{\lambda}$. There
 is a natural action of the basechange group $GL(\alpha) = \prod_{i=1}^k GL_{a_i}$ on this scheme
 and the corresponding quotient variety $\wis{iss}_{\alpha}~\Pi_{\lambda}$ represents the
 isomorphism classes of semisimple $\alpha$-dimensional representations of $\Pi_{\lambda}$. The
 main coadjoint orbit result of \cite{Ginzburg:2000} and \cite[Thm. 5.5]{BocklandtLB:2000} is
 
 \begin{theorem} \label{coadjoint} If $\alpha$ is a minimal element of $\Sigma_{\lambda}$, the set of dimension vectors
 of simple representations of $\Pi_{\lambda}$, then $\wis{iss}_{\alpha}~\Pi_{\lambda}$ is a
 coadjoint orbit for the necklace Lie algebra $\mathfrak{N}_Q$.
 \end{theorem}
 
 The first description of $\Sigma_{\lambda}$ is due to W. Crawley-Boevey \cite{Crawley:1999}. In
 \cite{LB:2001} the author gave an alternative characterization. In \cite{BocklandtLB:2000} we gave
 a conjectural ringtheoretical explanation for these coadjoint orbit results in terms of
 noncommutative notions of smoothness which we will recall in the next section. The main result of
 this paper is the following affirmative solution to this conjecture.
 
 \begin{theorem} \label{main} The following are equivalent
 \begin{enumerate}
 \item{$\alpha$ is a minimal element of $\Sigma_{\lambda}$.}
 \item{$\wis{iss}_{\alpha}~\Pi_{\lambda}$ is a coadjoint orbit for $\mathfrak{N}_Q$.}
 \item{$\wis{iss}_{\alpha}~\Pi_{\lambda}$ is a smooth variety.}
 \item{$\int_{\alpha}~\Pi_{\lambda}$ is an Azumaya algebra over the smooth variety $\wis{iss}_{\alpha}~\Pi_{\lambda}$.}
  \item{$\Pi_{\lambda}$ is $\alpha$-smooth in the sense of Procesi \cite{Procesi:1987}.}
 \end{enumerate}
 \end{theorem}
 
 The outline of this paper, as well as the proof of this result is summarized in the following picture
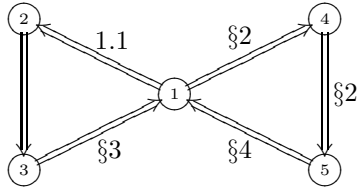
\[
\begin{xy}
\POS (0,0) *\cir<6pt>{}*+{\txt\tiny{$1$}} ="S1",
\POS (-20,10) *\cir<6pt>{}*+{\txt\tiny{$2$}} ="S2",
\POS (-20,-10) *\cir<6pt>{}*+{\txt\tiny{$3$}} ="S3",
\POS (20,10) *\cir<6pt>{}*+{\txt\tiny{$4$}} ="S4",
\POS (20,-10) *\cir<6pt>{}*+{\txt\tiny{$5$}} ="S5",
\POS"S1" \ar@{=>}_{\txt{1.1}} "S2"
\POS"S1" \ar@{=>}^{\txt{\S 2}} "S4"
\POS"S2" \ar@{=>} "S3"
\POS"S3" \ar@{=>}_{\txt{\S 3}} "S1"
\POS"S4" \ar@{=>}^{\txt{\S 2}} "S5"
\POS"S5" \ar@{=>}^{\txt{\S 4}} "S1"
\end{xy}
\]

\par \vskip 3mm

 \section{Noncommutative smoothness}
 
Path algebras of quivers are examples of formally smooth algebras as defined and studied in
\cite{CuntzQuillen:1995}, that is they have the lifting property for algebra morphisms 
modulo  nilpotent ideals. As a consequence they have a good theory of differential forms,
see for example \cite{BocklandtLB:2000}. If $A$ is a formally smooth $V$-algebra, then for each
dimension vector $\alpha$ the scheme of $\alpha$-dimensional representations $\wis{rep}_{\alpha}~A$
is smooth and noncommutative (relative) differential forms of $A$ induce ordinary $GL(\alpha)$-invariant
differential forms on these manifolds and hence on the corresponding quotient varieties
$\wis{iss}_{\alpha}~A$.

On the other hand we will see in the next sections that the deformed preprojective algebra
$\Pi_{\lambda}$ is {\em not} formally smooth as many of the representation schemes
$\wis{rep}_{\alpha}~\Pi_{\lambda}$ are singular. The quotient $\C \bar{Q} \rOnto \Pi_{\lambda}$
indicates that $\Pi_{\lambda}$ corresponds to a singular noncommutative subscheme
of the noncommutative manifold corresponding to the formally smooth algebra $\C \bar{Q}$. As a
consequence, the differential forms of $\C \bar{Q}$, when restricted to the singular subvariety
$\Pi_{\lambda}$, may have rather unpredictable behavior.

Still, it may be that the induced $GL(\alpha)$-invariant differential forms on some of the
representation schemes $\wis{rep}_{\alpha}~\Pi_{\lambda}$ have desirable properties, in particular
if $\wis{rep}_{\alpha}~\Pi_{\lambda}$ is a smooth variety. For this reason we need a notion
of noncommutative smoothness relative to a specific dimension vector $\alpha$. This notion was
introduced by C. Procesi in \cite{Procesi:1987} and investigated further in \cite{LB:2000}. We
briefly recall the definition and main results from \cite{Procesi:1987}.

With $\wis{alg@}\alpha$ we denote the category of $V$-algebras $C$ equipped with a $V$-linear trace map
$C \rTo^t C$ satisfying $t(ab) = t(ba)$, $t(a)b = bt(a)$, $t(t(a)b) = t(a)t(b)$ for all $a,b \in C$
and such that $t(e_i) = a_i$ if $\alpha = (a_1,\hdots,a_k)$ and $C$ satisfies the formal Cayley-Hamilton
identity of degree $n$ where $n = \sum_i a_i$. To explain the last definition, consider the
characteristic polynomial $\chi_M(t)$ of a general $n \times n$ matrix $M$ which is a polynomial
in a central variable $t$ with coefficients which can be expressed as polynomials with rational
coefficients in $Tr(M),Tr(M^2),\hdots,Tr(M^n)$. Replacing $M$ by $a$ and $Tr(M^i)$ by $t(a^i)$ we
have a formal characteristic polynomial $\chi_a(t) \in C[t]$ and we require that
$\chi_a(a) = 0$ for all $a \in C$. Morphisms in $\wis{alg@}\alpha$ are $V$-algebra morphisms which
are trace preserving.

An algebra $C$ in $\wis{alg@}\alpha$ is said to be $\alpha$-smooth if it satisfies
the lifting property for morphisms modulo nilpotent ideals in $\wis{alg@}\alpha$. That is, every diagram
\[
\begin{diagram}
B & \rOnto^{\pi} & \frac{B}{I} \\
& \luDotsto_{\exists \tilde{\phi}} & \uTo^{\phi} \\
& & C
\end{diagram}
\]
with $B,\tfrac{B}{I}$ in $\wis{alg@}\alpha$, $I$ a nilpotent
ideal and $\pi$ and $\phi$ trace preserving
maps, can be completed with a trace preserving algebra map $\tilde{\phi}$.

The forgetful functor $\wis{alg@}\alpha \rTo \wis{alg}$ has a left inverse $\int_{\alpha}$ which
assigns to a $V$-algebra $A$ the algebra $\int_{\alpha}~A$ obtained by formally adjoining traces
to $A$ and then modding out all Cayley-Hamilton identities of degree $n$. From \cite{Procesi:1987}
we recall geometric reconstruction results for $\int_{\alpha}~A$ and its central subalgebra
$\oint_{\alpha}~A = t~\int_{\alpha}~A$ as well as the characterization of $\alpha$-smoothness.

\begin{theorem}[C. Procesi] \label{procesi} With notations as above we have
\begin{enumerate}
\item{The algebra $\int_{\alpha}~A$ is the ring of $GL(\alpha)$-equivariant maps from
$\wis{rep}_{\alpha}~A$ to $M_n(\C)$ where $GL(\alpha)$ acts on the latter by conjugation via the
diagonal embedding $GL(\alpha) \rInto GL_n$, that is,
\[
\int_{\alpha}~A = M_n(\C[\wis{rep}_{\alpha}~A])^{GL(\alpha)} \]}
\item{The image of the trace map on $\int_{\alpha}~A$ is the ring of $GL(\alpha)$-invariant
polynomial functions on $\wis{rep}_{\alpha}~A$, that is
\[
\oint_{\alpha}~A = \C[\wis{iss}_{\alpha}~A] \]}
\item{$A$ is $\alpha$-smooth in the sense of Procesi, that is, $\int_{\alpha}~A$ is $\alpha$-smooth
in $\wis{alg@}\alpha$
if and only if $\wis{rep}_{\alpha}~A$ is a smooth variety.}
\end{enumerate}
\end{theorem}
 
Recall that an algebra $C$ in $\wis{alg@}\alpha$ is said to be an Azumaya algebra if and only if
every trace preserving morphism $C \rTo M_n(\C)$ is an epimorphism. If we start with a $V$-algebra
$A$, then a trace preserving algebra map $\int_{\alpha}~A \rTo M_n(\C)$ corresponds one-to-one
to a $V$-algebra map $A \rTo M_n(\C)$ hence to a geometric point of $\wis{rep}_{\alpha}~A$.
The Azumaya property for $\int_{\alpha}~A$ is therefore equivalent to saying that the
quotient map
\[
\wis{rep}_{\alpha}~A \rOnto^{\pi} \wis{iss}_{\alpha}~A \]
is a principal $PGL(\alpha)$-fibration in the \'etale topology. For, in general a geometric point in
$\wis{iss}_{\alpha}~A$ determines an isomorphism class of a semi-simple $\alpha$-dimensional representation
of $A$ and the map
$\pi$ sends an $\alpha$-dimensional representation to the direct sum of its Jordan-H\"older
factors. 

\begin{proposition} The following implications of theorem~\ref{coadjoint} hold :

$(1) \Rightarrow (4)$ : If $\alpha$ is a minimal element of $\Sigma_{\lambda}$, then $\int_{\alpha}~\Pi_{\lambda}$ is
an Azumaya algebra over the smooth variety $\wis{iss}_{\alpha}~\Pi_{\lambda}$.

$(4) \Rightarrow (5)$ : If $\int_{\alpha}~\Pi_{\lambda}$ is an Azumaya algebra over the smooth variety
$\wis{iss}_{\alpha}~\Pi_{\lambda}$, then $\Pi_{\lambda}$ is $\alpha$-smooth in the sense of Procesi.
\end{proposition}

\begin{proof}
$(1) \Rightarrow (4)$ : Consider the {\em complex moment map}
\[
\wis{rep}_{\alpha}~\bar{Q} \rTo^{\mu_{\C}} M_{\alpha}^0(\C) \qquad V \mapsto \sum_{a
\in Q_a} [V_a,V_{a^*}]
\]
where $M_{\alpha}^0(\C)$ is the subspace of $k$-tuples $(m_1,\hdots,m_k)
\in M_{a_1}(\C) \oplus \hdots
\oplus M_{a_k}(\C)$ such that $\sum_i tr(m_i) = 0$.
For $\lambda
= (\lambda_1,\hdots,\lambda_k) \in \Q^k$
such that $\sum_i a_i \lambda_i = 0$ we consider the element
$\underline{\lambda} = (\lambda_1 \een_{n_1},\hdots,\lambda_k \een_{n_k})$
in $M_{\alpha}^0(\C)$. The inverse image
$\mu^{-1}_{\C}(\underline{\lambda}) =
\wis{rep}_{\alpha}~\Pi_{\lambda}$.
By a result of M. Artin \cite{Artin:1969} one knows that the geometric
points of the quotient scheme
$\wis{iss}_{\alpha}~\Pi_{\lambda}$ are the isomorphism
classes of $\alpha$-dimensional
semi-simple representations of $\Pi_{\lambda}$. Because $\alpha$ is a minimal
element of $\Sigma_{\lambda}$ all $\alpha$-dimensional representations of $\Pi_{\lambda}$
must be simple (consider the dimension vectors of Jordan-H\"older components) so each
fiber of the quotient map $\pi$ is isomorphic to $PGL(\alpha)$. The fact that $\mu^{-1}_{\C}(\underline{\lambda})$ is smooth
if $\alpha$ is a minimal non-zero element of $\Sigma_{\lambda}$ follows from computing the
differential of the complex moment map, see also \cite[lemma 5.5]{Crawley:1999}. Because
$\wis{rep}_{\alpha}~\Pi_{\lambda} \rOnto^{\pi} \wis{iss}_{\alpha}~\Pi_{\lambda}$ is a principal
$PGL(\alpha)$-fibration, $\int_{\alpha}~A$ is an Azumaya algebra and as
the total space $\wis{rep}_{\alpha}~\Pi_{\lambda}$ is smooth it follows that also the
base space $\wis{iss}_{\alpha}~\Pi_{\lambda}$ is smooth. 

$(4) \Rightarrow (5)$ : If $\int_{\alpha}~\Pi_{\lambda}$ is an Azumaya algebra, it follows that
\[
\wis{rep}_{\alpha}~\Pi_{\lambda} \rOnto^{\pi} \wis{iss}_{\alpha}~\Pi_{\lambda} \]
is a principal $PGL(\alpha)$-fibration. If in addition the basespace is smooth, so is the top
space $\wis{rep}_{\alpha}~\Pi_{\lambda}$. The assertion follows from Procesi's characterization
of $\alpha$-smoothness, theorem~\ref{procesi}.
\end{proof}

\section{Central singularities}

Clearly, if $\wis{iss}_{\alpha}~\Pi_{\lambda}$ is a coadjoint orbit of $\mathfrak{N}_Q$ it is
a smooth variety. In this section we will show that unless $\alpha$ is a minimal element of
$\Sigma_{\lambda}$ the quotient variety $\wis{iss}_{\alpha}~\Pi_{\lambda}$ always has
singularities. The crucial ingredient in the proof is the \'etale local description of
$\wis{iss}_{\alpha}~\Pi_{\lambda}$ due to W. Crawley-Boevey \cite{Crawley:2001}.

Let $\chi_Q$ be the Euler-form of the quiver $\vec{Q}$, that is, the bilinear form
\[
\chi_Q~:~\Z \times \Z \rTo \Z \]
is determined by the matrix $(\chi_{ij})_{i,j} \in M_k(\Z)$ where $\chi_{ij} = \delta_{ij} -
\# \{ a \in Q_a~\text{starting at $v_i$ and ending in $v_j$}~\}$. We denote the
symmetrization of $\chi_Q$ by $T_Q$ (the Tits form) and $p(\alpha) = 1 - \chi_Q(\alpha,\alpha)$
for every dimension vector $\alpha$.

Throughout we assume that $\alpha \in \Sigma_{\lambda}$ and we consider a geometric point
$\xi \in \wis{iss}_{\alpha}~\Pi_{\lambda}$ which determines an isomorphism class of a
semisimple $\alpha$-dimensional representation of $\Pi_{\lambda}$ say
\[
M_{\xi} = S_1^{\oplus e_1} \oplus \hdots \oplus S_u^{\oplus e_u} \]
where the $S_i$ are distinct simple representations of $\Pi_{\lambda}$ with dimension vectors
$\beta_i$ and occurring in $M_{\xi}$ with multiplicity $e_i$, that is,
$\alpha = \sum_{i=1}^u e_i \beta_i$. We say that $\xi$ is of representation type
$\tau = (e_1,\beta_1;\hdots;e_u,\beta_u)$.

Construct a new (symmetric) quiver $\Gamma_{\tau}$ on $u$ vertices $\{ v'_1,\hdots,v'_u \}$
(corresponding to the distinct simple components) such that there are
\begin{itemize}
\item{$2p(\beta_i)$ loops in vertex $v'_i$, and }
\item{$-T_Q(\beta_i,\beta_j)$ directed arrows from $v'_i$ to $v'_j$.}
\end{itemize}
We also consider the dimension vector $\alpha_{\tau} = (e_1,\hdots,e_u)$ for $\Gamma_{\tau}$.

\begin{theorem}[W. Crawley-Boevey] \label{crawley} With notations as above there is an \'etale isomorphism between
\begin{enumerate}
\item{a neighborhood of $\xi$ in $\wis{iss}_{\alpha}~\Pi_{\lambda}(\bar{Q})$, and}
\item{a neighborhood of the trivial representation $\overline{0}$ in $\wis{iss}_{\alpha_{\tau}}~\Pi_0(\Gamma_{\tau})$}
\end{enumerate}
where $\Pi_0(\Gamma_{\tau})$ is the preprojective algebra corresponding to the double quiver
$\Gamma_{\tau}$.
\end{theorem}

In particular it follows that $\wis{iss}_{\alpha}~\Pi_{\lambda}$ is smooth in all points $\xi$ 
of representation type $\tau = (1,\alpha)$ (the so called {\em Azumaya locus}) and that the dimension of $\wis{iss}_{\alpha}~\Pi_{\lambda}$
is equal to $2p(\alpha)$. In \cite[Prop 8.6]{Crawley:2001} it was proved that $\wis{iss}_{\alpha}~\Pi_0(\bar{Q})$ has
singularities in case $\vec{Q}$ is a quiver without loops and $\alpha$ is an imaginary indivisible
root of $\Sigma_0$. Recall that a dimension vector is said to be indivisible if the greatest common divisor of
its components is one. 

\begin{theorem} For $\alpha \in \Sigma_{\lambda}$, the smooth locus of $\wis{iss}_{\alpha}~\Pi_{\lambda}$
coincides with the Azumaya locus.

In particular, if $\alpha$ is not a minimal element of $\Sigma_{\lambda}$, then $\wis{iss}_{\alpha}~\Pi_{\lambda}$
is singular, that is, implication $(3) \Rightarrow (1)$ of theorem~\ref{main} holds.
\end{theorem}

\begin{proof} With $\wis{iss}_{\alpha}(\tau)$ we will denote the locally closed subvariety of
$\wis{iss}_{\alpha}~\Pi_{\lambda}$ consisting of all geometric points $\xi$ of representation type
$\tau$. Observe that there is a natural ordering on the set of representation types 
\[
\tau \leq \tau'~\quad~\Longleftrightarrow~\quad \wis{iss}_{\alpha}(\tau) \subset
\overline{\wis{iss}_{\alpha}(\tau')} \]
where the closure is with respect to the Zariski topology. Clearly, if we can prove that all
points of $\wis{iss}_{\alpha}(\tau')$ are singular then so are those of $\wis{iss}_{\alpha}(\tau)$.

Let $\xi$ be a point outside of the Azumaya locus of representation type $\tau = (e_1,\beta_1;\hdots;
e_u,\beta_u)$ then by theorem~\ref{crawley} is suffices to prove that $\wis{iss}_{\alpha_{\tau}}~\Pi_0(\Gamma_{\tau})$
is singular in $\overline{0}$.

Assume that $\Gamma_{\tau}$ has $2p(\beta_i) > 0$ loops in the vertex $v'_i$ where $e_i > 1$. This
means that there are infinitely many nonisomorphic simple $\Pi_{\lambda}$-representations of
dimension vector $\beta_i$, but then $\tau < \tau'$ where
\[
\tau' = (e_1,\beta_1;\hdots;e_{i-1},\beta_{i-1};\underbrace{1,\beta_i;\hdots;1,\beta_i}_{e_i};e_{i+1},\beta_{i+1};\hdots;e_u,\beta_u) \]
and by the above remark it suffices to prove singularity for $\tau'$. That is, we may assume that the
quiver setting $(\Gamma_{\tau},\alpha_{\tau})$ is such that the symmetric quiver $\Gamma_{\tau}$
has loops only at vertices $v'_i$ where the dimension $e_i = 1$.

Assume moreover that $\alpha_{\tau}$ is indivisible (which by the above can be arranged
once we start from
a type $\tau$ such that $\Gamma_{\tau}$ has loops). Recall from \cite{LeBruynProcesi:1990}
that invariants of quivers are generated by traces along oriented cycles in the quiver. As
a consequence we have algebra generators of the coordinate ring $\C[\wis{iss}_{\alpha_{\tau}}~\Pi_0(\Gamma_{\tau})]$
which is a graded algebra by homogeneity of the defining relations of the preprojective
algebra $\Pi_0(\Gamma_{\tau})$. To prove singularity it therefore suffices to prove that the
coordinate ring is not a polynomial ring. Let $\Gamma_{\tau}'$ be the quiver obtained from
$\Gamma_{\tau}$ by removing all loops (which by the above reduction exist only at vertices where the
dimension is one). Because the relations of the preprojective algebra are irrelevant for loops in such
vertices we have that $\C[\wis{iss}_{\alpha_{\tau}}~\Pi_0(\Gamma_{\tau})]$ is a polynomial ring
(in the variables corresponding to the loops) over $\C[\wis{iss}_{\alpha_{\tau}}~\Pi_0(\Gamma_{\tau}')]$.
By \cite[Prop 8.6]{Crawley:2001} we know that $\wis{iss}_{\alpha_{\tau}}~\Pi_0(\Gamma'_{\tau})$ is
singular, finishing the proof in this case.

The remaining case is when $\Gamma_{\tau}$ contains no loops (that is, all $\beta_i$ are real roots
for $\vec{Q}$) and when $\alpha_{\tau}$ is divisible. Because $\alpha_{\tau}$ is the dimension
vector of a simple representation of $\Pi_0(\Gamma_{\tau})$ we know from \cite{LB:2001} that the
quiver setting $(\Gamma_{\tau},\alpha_{\tau})$ is such that $\Gamma_{\tau}$ contains a subquiver
say on the vertices $T=\{ v'_{i_1},\hdots,v'_{i_z} \}$ which is the double of a tame quiver such that
$(\alpha_{\tau} \mid T) \geq \delta$ where $\delta = (d_{i_1},\hdots,d_{i_z})$ 
is the imaginary root of this tame subquiver.
Consider the representation type of $\alpha_{\tau}$ for $\Pi_0(\Gamma_{\tau})$
\[
\gamma = (1,\delta;e_1,\epsilon_1;\hdots;e_{i_1}-d_{i_1},\epsilon_{i_1},\hdots;
e_{i_x}-d_{i_x},\epsilon_{i_x};\hdots;e_u,\epsilon_u) \]
If we can show that a points is the $\gamma$-stratum of $\wis{iss}_{\alpha_{\tau}}~\Pi_0(\Gamma_{\tau})$
is singular, then the quotient scheme is singular in the trivial representation and we are done.
Consider the quiver $\Gamma_{\gamma}$, then it has loops in the vertex corresponding to $(1,\delta)$.
Moreover, $\alpha_{\gamma}$ is indivisible so we can repeat the argument above.
The fact that $\alpha_{\tau}$ was assumed
to be divisible asserts that $\gamma$ is not the Azumaya type $(1,\alpha_{\tau})$, finishing the
proof.
\end{proof}

\section{Noncommutative singularities}

We can refine the notion of $\alpha$-smoothness to allow for noncommutative singularities
with respect to the dimension vector $\alpha$. Let $A$ be a $V$-algebra and consider the
quotient map
\[
\wis{rep}_{\alpha}~A \rOnto^{\pi} \wis{iss}_{\alpha}~A \]
and consider the open subvariety $\wis{sm}_{\alpha}~A$ of $\wis{iss}_{\alpha}~A$
consisting of those geometric points $\xi$ such that $\wis{rep}_{\alpha}~A$ is smooth along the
fiber $\pi^{-1}(\xi)$ and call $\wis{sm}_{\alpha}~A$ the $\alpha$-smooth locus of $A$. In particular,
$A$ is $\alpha$-smooth in the sense of Procesi if and only if $\wis{sm}_{\alpha}~A = \wis{iss}_{\alpha}~A$.

Returning to $\Pi_{\lambda}$ it is clear from the foregoing that the Azumaya locus is contained in
the $\alpha$-smooth locus for $\alpha \in \Sigma_{\lambda}$. We will show in this section that this
these loci are actually identical showing that deformed preprojective algebras are as singular as
possible. This result should be compared to a similar result on quantum groups at roots
of unity \cite{LeBruyn:2000q}.

Using the notations needed in theorem~\ref{crawley} we observe that the method of proof actually
proves a stronger result which is a symplectic version of Luna slices \cite{Luna:1973}, see also
\cite[\S 41]{GS:1990}.

\begin{theorem} \label{etale} There is a $GL(\alpha)$-equivariant \'etale isomorphism between
\begin{enumerate}
\item{a neighborhood of the orbit of $M_{\xi}$ in $\wis{rep}_{\alpha}~\Pi_{\lambda}$, and}
\item{a neighborhood of the orbit of $\overline{(\een_{\alpha},0)}$ in the principal fiber bundle
\[
GL(\alpha) \times^{GL(\alpha_{\tau})} \wis{rep}_{\alpha_{\tau}}~\Pi_0(\Gamma_{\tau}) \]}
\end{enumerate}
\end{theorem}

We are now in a position to prove the final implication of theorem~\ref{main}.

\begin{theorem} If $\alpha \in \Sigma_{\lambda}$, then the $\alpha$-smooth locus
of $\Pi_{\lambda}$ coincides with the Azumaya locus.

In particular, if $\Pi_{\lambda}$ is $\alpha$-smooth in the sense of Procesi, then there
is only one representation type $(1,\alpha)$, that is, $\alpha$ is a minimal element of
$\Sigma_{\lambda}$. Hence, $(5) \Rightarrow (1)$ of theorem~\ref{main} holds.
\end{theorem}

\begin{proof} Assume that $\xi \in \wis{sm}_{\alpha}~\Pi_{\lambda}$ and of representation type
$\tau = (e_1,\beta_1;\hdots;e_u,\beta_u)$. Then, $\wis{rep}_{\alpha}~\Pi_{\lambda}$ is smooth
in a neighborhood of the {\em closed} $GL(\alpha)$-orbit of the semisimple representation
$M_{\xi}$ (closedness follows from \cite{Artin:1969}). By a result of Voigt \cite{Kraftbook}
we know that the normalspace $N_{\xi}$ to the orbit in $\wis{rep}_{\alpha}~\Pi_{\lambda}$ is
equal to the space of self-extensions
\[
Ext^1_{\Pi_{\lambda}}(M_{\xi},M_{\xi}) = \oplus_{i,j=1}^u Ext^1_{\Pi_{\lambda}}(S_i,S_j)^{\oplus e_ie_j} \]
As a consequence we can identify this space with the representation space $\wis{rep}_{\alpha_{\tau}}~\Delta_{\tau}$
where $\Delta_{\tau}$ is the quiver on $u$ vertices $\{ w'_1,\hdots,w'_u \}$ having
\begin{itemize}
\item{$dim~Ext^1_{\Pi_{\lambda}}(S_i,S_i)$ loops in vertex $w_i'$, and}
\item{$dim~Ext^1_{\Pi_{\lambda}}(S_i,S_j)$ directed arrows from $w_i'$ to $w_j'$}
\end{itemize}
Moreover, the action of the stabilizer subgroup of $M_{\xi}$ (which is $GL(\alpha_{\tau})$) on the
normal; space is the
basechange action of this group on $\wis{rep}_{\alpha_{\tau}}~\Delta_{\tau}$.
By the Luna slice theorem \cite{Luna:1973} we have a $GL(\alpha)$-equivariant \'etale isomorphism
between
\begin{enumerate}
\item{a neighborhood of the orbit of $M_{\xi}$ in $\wis{rep}_{\alpha}~\Pi_{\lambda}$, and}
\item{a neighborhood of the orbit of $\overline{(\een_{\alpha},0)}$ in the principal fiber bundle
\[
GL(\alpha) \times^{GL(\alpha_{\tau})} \wis{rep}_{\alpha_{\tau}}~\Delta_{\tau} \]}
\end{enumerate}
Combining this \'etale description with the one from theorem~\ref{etale} we deduce an \'etale
$GL(\alpha_{\tau})$-isomorphism between the representation scheme $\wis{rep}_{\alpha_{\tau}} \Pi_0(\Gamma_{\tau})$
(in a neighborhood of the trivial representation) and the representation space $\wis{rep}_{\alpha_{\tau}}~\Delta_{\tau}$
(in a neighborhood of the trivial representation). 

But then $\overline{0} \in \wis{sm}_{\alpha_{\tau}}~\Pi_0(\Gamma_{\tau})$ and in \cite[Thm. 6.3]{BocklandtLB:2000}
it was shown that for a preprojective algebra the smooth locus coincides with the Azumaya algebra.
The only way the trivial representation can be a simple representation of $\Pi_0(\Gamma_{\tau})$
(or indeed, even of $\C \Gamma_{\tau}$) is when $\Gamma_{\tau}$ has only one vertex and the
dimension vector is $\alpha_{\tau} = 1$. But then the representation type of $\xi$ is $(1,\alpha)$
finishing the proof.
\end{proof}


\begin{thebibliography}{19}

\bibitem{Artin:1969}~M. Artin, {\it On Azumaya algebras and finite 
dimensional representations
of rings}, J.Alg. 11 (1969) 523-563


\bibitem{BerestWilson:2000}~Yu. Berest and G. Wilson, {\it 
Automorphisms and ideals of the Weyl
algebra} preprint, London (1999) see also {\tt math.QA/0102190} (2001)

\bibitem{BocklandtLB:2000}~R. Bocklandt and L. Le Bruyn, {\it Necklace Lie algebras and
noncommutative symplectic geometry} {\tt math.AG/0010030} (2000), Math. Z. (to appear)


\bibitem{Crawley:1999}~W. Crawley-Boevey {\it Geometry of the moment 
map for representations
of quivers} Compositio Math. 126 (2001) 257-293

\bibitem{Crawley:2001}~W. Crawley-Boevey {\it Normality of Marsden-Weinstein reductions
for representations of quivers} {\tt math.AG/0105247} (2001)


\bibitem{CuntzQuillen:1995}~J. Cuntz and D. Quillen {\it Algebra 
extensions and nonsingularity}
Journal AMS 8 (1995) 251-289

\bibitem{GS:1990}~V. Guillemin and S. Sternberg {\it Symplectic techniques in physics}
Cambridge University Press (1990)

\bibitem{Ginzburg:1999}~V. Ginzburg {\it Non-commutative symplectic 
geometry and Calogero-Moser
space} preprint Chicago, preliminary version (1999)

\bibitem{Ginzburg:2000}~V. Ginzburg {\it Non-commutative symplectic 
geometry, quiver varieties
and operads} preprint Chicago (2000) {\tt math.QA/0005165} (2000)


\bibitem{Kontsevich:1993}~M. Kontsevich {\it Formal non-commutative 
symplectic geometry}
Gelfand seminar 1990-1992, Birkhauser (1993) 173-187

\bibitem{Kraftbook}~H.P. Kraft {\it Geometrische Methoden in der Invariantentheorie}
Aspekte der mathematik, vol D1 (1984) Vieweg-Verlag, Braunschweig

\bibitem{LeBruyn:2000q}~L. Le Bruyn {\it The singularities of quantum groups}
Proc. London Math. Soc. 80 (2000) 304-336

\bibitem{LB:2000}~L. Le Bruyn {\it Local structure of Schelter-Procesi smooth orders},
Trans. Amer. Math. Soc. 352 (2000) 4815-4841

\bibitem{LB:2001}~L. Le Bruyn {\it Simple roots of deformed preprojective algebras}
{\tt math.RA/0107027} (2001)

\bibitem{LeBruynProcesi:1990}~L. Le Bruyn and C. Procesi, {\it Semisimple representations of quivers}
Trans. AMS 317 (1990) 585-598

\bibitem{Luna:1973}~D. Luna, {\it Slices etales} Bull.Soc.Math. France Mem 33 
(1973) 81-105


\bibitem{Procesi:1987}~C. Procesi {\it A formal inverse to the 
Cayley-Hamilton theorem}
J.Alg. 107 (1987) 63-74


\end{thebibliography}
\end{document}